\documentclass[11pt]{article}
\usepackage{amsmath,amsthm,amssymb,xypic,graphicx}
\input xy
\xyoption{all}
\newtheorem{definition}{Definition}
\newtheorem{theorem}{Theorem}
\newtheorem{lemma}{Lemma}
\newtheorem{conjecture}{Conjecture}
\newtheorem{proposition}{Proposition}

\newcommand{\frakm}{\mathfrak{m}}
\newcommand{\frakp}{\mathfrak{p}}

\newcommand{\calH}{\mathcal H}
\newcommand{\calL}{\mathcal L}
\newcommand{\calU}{\mathcal U}
\newcommand{\calV}{\mathcal V}
\newcommand{\hB}{B}
\newcommand{\del}{\ensuremath{\partial}}

\newcommand{\harm}{\ensuremath{\mathcal{H}_{(2)}}}

\newcommand{\CC}{\mathbb C}

\newcommand{\CY}{C_{\phi}Y}
\newcommand{\calD}{\mathcal D}
\newcommand{\ran}{{\mathrm{ran}}}

\def\slashs{\,\hbox{\slash}\kern-6.5pt S}
\def\slashd{\,\hbox{\slash}\kern-8.0pt D}

\begin{document}
\title{Hodge and signature theorems for a family of manifolds with fibration boundary}
\author{Eugenie Hunsicker \thanks{Partially supported by 
the NSF through an ROA supplement to grant DMS-0204730} 
\\ Lawrence University}
  
\maketitle

\section{Introduction}

The Hodge theorem and Hirzebruch signature theorem form an 
important bridge between geometric and topological properties of 
compact smooth manifolds.  There has been a great deal of work over the 
past thirty years aimed at understanding how to generalize these theorems
to $L^2$ results in the noncompact and singular settings.  
Early and important work was done by Atiyah, Patodi and Singer \cite{APS1}.  
Their work concerned both manifolds with the simplest sort of singularities, 
namely boundaries, and noncompact manifolds with cylindrical ends, 
that is, manifolds which off a compact set are isometric to $(0, \infty) \times N$ 
for some compact manifold $N$.  They proved both a Hodge result and a signature result.  
Their Hodge result says that the space of $L^2$ harmonic forms on a 
manifold, $\hat{M}$, with cylindrical end is canonically isomorphic to the image 
of relative cohomology of  $\hat{M}$ in its absolute cohomology, i.e.
\begin{equation}
\harm^k(\hat{M},g) \cong \mbox{Im}(H^k(\hat{M}, N) \longrightarrow H^k(\hat{M})).
\label{eq:apshodge}
\end{equation}
Their signature result says that the topological signature of a manifold with 
boundary, $\overline{M}$, is equal first of all to the signature of the pairing 
on middle degree $L^2$ harmonic forms given by integration, and secondly that both satisfy
\begin{equation}
\sigma(M) = \int_M L(p(M, N)) - \eta(N)(0),
\label{eq:apssign}
\end{equation}
where $L$ is the Hirzebruch $L$-polynomial in the relative Pontrjagin classes on $M$ and
$\eta$ is a spectral invariant of the boundary, $N$, of $M$. 

The connection between $L^2$ Hodge theorems and intersection cohomology 
was made by Cheeger shortly after the dual intersection homology theoreies 
had been defined by Goresky and MacPherson \cite{GM1}.  Cheeger showed 
in  \cite{Ch} that for a pseudomanifold, $X$, with conical singularities and only even 
codimensional strata, the space of $L^2$ harmonic forms on the regular set is isomorphic 
to the middle perversity intersection cohomology, the unique intersection cohomology 
that satisfies Poincar\'{e} duality, that is:
\begin{equation}
H^i_2(X^{reg}, g_c) \cong IH^i_{m}(X),
\label{eq:cheegerhodge}
\end{equation} 

In this paper, Cheeger also noted the importance in $L^2$ Hodge theorems 
for incomplete manifolds generally of considering 
different "ideal boundary conditions".  That is, whereas in the complete case, there is a unique
closed extension of the exterior derivative, $d$, in the incomplete case, there may be 
several different closed extensions.  These correspond to several different $L^2$ cohomologies and
several different closed self-adjoint extensions of the Laplace operator.  
For instance, in the case of a manifold with boundary, $M$, the exterior derivative may have a number
of closed extensions interpolating between the so called maximal and minimal extensions.  The
cohomology corresponding to the maximal extension of $d$ is absolute 
cohomology on $M$.  Its classes are naturally represented by $L^2$ harmonic forms satisfying Neumann boundary conditions.   
The complex given by the minimal extension of $d$ generates relative cohomology on $M$.   Its
classes are naturally represented by harmonic $L^2$ 
forms satisfying Dirichlet boundary conditions.  The even codimension condition 
on conical pseudomanifolds in Cheeger's paper avoids this complication.

In the 1990's, new work was done on the eta invariant defined
by Atiyah, Patodi and Singer.  A number of mathematicians began to study its behavior
in a family of fibre bundle metrics which become singular.  Bismut and Cheeger, [5], 
Melrose and Mazzeo, [38], and Dai, [14], all studied the eta invariant under such
adiabatic limits. In 1990, Muller linked his work on signature theorems for manifolds with cusps of
rank 1 to this new work on the eta invariant, [49], and interpreted it in terms of middle
perversity intersection cohomology.  He showed that the $L^2$ signature for such a manifold
was equal to the $L^2$ signature for the manifold with boundary formed by cutting off the cusps
plus the $\tau$ invariant of the resulting boundary fibration, defined by Dai.

Recently, Dai and Cheeger followed up Cheeger's Hodge theorem with
a signature theorem for conical pseudomanifolds with one even codimensional
smooth singular stratum. They show in \cite{CD} that
the signature of the intersection pairing on middle degree $L^2$ forms for a 
conical pseudomanifold, $X$, whose singular set, $B$, is a smooth compact even codimensional
manifold is given by the formula:
\begin{equation}
\sigma(X)= \mbox{sgn}\,
\mbox{Im}\,\big( H^*(X, B) \longrightarrow H^*(X^{reg}))\big) +  \tau.
\label{eq:conesign}
\end{equation}
where $\tau$ is an invariant of the fibre bundle over $B$ which forms the boundary of the normal
neighborhood of the singular stratum.   Recently, in \cite{HM}, the author and Mazzeo have extended Cheeger and Dai's results to conical pseudomanifolds with one singular stratum in the form of a smooth compact manifold which may be odd dimensional.  

The similarity and relationships between the $L^2$ Hodge and signature results for manifolds
with boundary and manifolds with conical singularities leads to to the question of what happens as
we interpolate between these metrics on a given incomplete smooth manifold.  In this paper, we
generalize the techniques of \cite{CD} to answer this question.

Consider a manifold $\overline{M}$ with boundary $\del \overline{M} = Y$ where 
$Y \stackrel{\phi}{\rightarrow} B$ is a fibre bundle with fibre $F$.  Endow 
$M=\overline{M} - \del \overline{M}$ with a metric $g_c$ which is quasi-isometric near 
the boundary to one of the form
\begin{equation}
ds_c^2 = dr^2 + r^{2c}\tilde{h} + \phi^*ds_B^2
\label{eq:metrics}
\end{equation}
\noindent
where $\tilde{h}$ is a two form which restricts to a metric on each fibre of $Y$ and $0 \leq c \leq 1$.
Note that these metrics interpolate between a cylindrical metric when $c=0$ and a cone bundle metric
when $c=1$.  There is a second natural compactification, $X$, for the manifold $M$, obtained by 
collapsing the fibres of the boundary fibration and adding a compactifying copy of the base.

As mentioned before, on a complete manifold, the exterior derivative $d$ on compactly 
supported smooth forms and its formal 
adjoint $\delta$ each have a unique extension to $L^2$ forms, and the Hodge Laplacian given by these
extensions, $\overline{d}+\overline{\delta}$, is self-adjoint.  On an incomplete manifold, this is not 
generally true.  There may be many closed extensions of $d$, and many self-adjoint extensions of the 
Hodge Laplacian.  Perhaps the most natural of these are the minimal and maximal extensions of $d$ 
and $\delta$ and their associated self-adjoint Laplacians.  Recall that the minimal and maximal 
extensions of the exterior derivative, $d$, on compact smooth forms to $L^2$ forms on $(M, g_c)$ are 
defined by:

\begin{definition}  Let $\sigma \in L^2\Omega^k(M,g)$ for 
any $k$.  If there exists an 
$\eta \in L^2\Omega^{k+1}(M,g)$ such that for all $\phi 
\in C^{\infty}_0\Omega^{k+1}(M,g)$ we have
\begin{equation}
\langle \sigma , \delta \phi \rangle_{L^2(M)} = \langle 
\eta , \phi \rangle_{L^2(M)},
\label{eq:dmax}
\end{equation}
we say that $\sigma \in \mathcal{D}(d_{max})$ and we write 
$d_{max} \sigma = \eta$.
\end{definition}

\begin{definition} 
If $\sigma \in L^2\Omega^{k}(M,g)$ for any $k$, and if 
there is a sequence $\{\sigma_n\} \subset 
C^{\infty}_0\Omega^{k}(M,g)$ such that $\sigma = 
\lim_{L^2} \sigma_n$ and $\zeta = \lim_{L^2} d
\sigma_n$, we say $\sigma \in \mathcal{D}(d_{min})$ and 
$\zeta = d_{min}\sigma$.  In this case, it is
always true that $\sigma \in \mathcal{D}(d_{max})$ as 
well, and $d_{max} \sigma = d_{min} \sigma$.
\end{definition}

The minimal and maximal extensions of $\delta$ are defined analogously.  Both
extensions of $d$ define cohomology groups on $M$, which we will denote by $H^*_{min}(M, g_c)$
and $H^*_{max}(M, g_c)$.  We define two related self-adjoint extensions of the Hodge Laplacian 
$D= d+\delta$:  
\begin{equation}
D_{m,M}=d_{min}+\delta_{max}
\label{eq:mMlaplacian}
\end{equation}
\begin{equation}
D_{M,m}=d_{max}+\delta_{min}.
\label{eq:Mmlaplacian}
\end{equation}
We will 
denote the corresponding spaces of harmonic forms by $\mathcal{H}^*_{m,M}(M, g_c)$ and 
$\mathcal{H}^*_{M,m}(M,g_c)$, respectively.  

Our first theorem is a Hodge theorem for these cohomology groups:

\begin{theorem}
Let $\overline{M}$ be a manifold with boundary $\del \overline{M} = Y$ where 
$Y \stackrel{\phi}{\rightarrow} B$ is a fibre bundle with $f$-dimensional fibre $F$.  Endow 
$M=\overline{M} - \del \overline{M}$ with a metric $g_c$ which is quasi-isometric near 
the boundary to one of the form
$$
ds_c^2 = dr^2 + r^{2c}\tilde{h} + \phi^*ds_B^2
$$
\noindent
where $h$ is a two form which restricts to a metric on each fibre of $Y$.
Let $X$ be the compactification of $M$ obtained by collapsing
the fibres of the boundary fibration.  Then 
$$
H^*_{min}(M, g_c) \cong  \mathcal{H}^*_{m,M}(M, g_c) \cong 
\left\{ \begin{array}{ll}
I\!H_{\underline{\frakm}+[[1+ \frac{1}{2c}]]}^*(X, B) & f \mbox{ is even} \\
I\!H_{\underline{\frakm}+[[\frac{1}{2} + \frac{1}{2c}]]}^*(X, B) & f \mbox{ is odd} \\
\end{array} \right. 
$$
and
$$
H^*_{max}(M, g_c) \cong  \mathcal{H}^*_{M, m}(M, g_c) \cong 
\left\{ \begin{array}{ll}
I\!H_{\overline{\frakm}-[[1+ \frac{1}{2c}]]}^*(X, B) & f \mbox{ is even} \\
I\!H_{\overline{\frakm}-[[\frac{1}{2} + \frac{1}{2c}]]}^*(X, B) & f \mbox{ is odd} \\
\end{array} \right. ,
$$
\noindent
where $[[x]]$ denotes the greatest integer strictly less than $x$ and $\underline{\frakm}$ and 
$\overline{\frakm}$ are the two middle perversities for $X$. 
\label{th:hodge}
\end{theorem}
Here we use the notation $I\!H_{\mathfrak{p}}^*(X, B)$ instead of  $I\!H_{\mathfrak{p}}^*(X)$ in order to 
indicate a slightly more general definition of intersection cohomology than is 
standard.  In particular, it allows us to include the case where the boundary fibration
fibre $F$ is trivial, so $X$ is our original manifold with boundary, $\overline{M}$.

Note that when $c=1$, we get $[[\frac{1}{2} + \frac{1}{2c}]] =0$, so this result reduces
to the result for manifolds with edges in \cite{HM}, and in the case that the fibre is even dimensional,
it reduces to Cheeger's result in \cite{Ch}.  In the extended definition
of intersection cohomology , for $c$ sufficiently close to $0$, 
these spaces become relative and absolute cohomologies of $M$, respectively, thus reducing to
the known results for manifolds with boundary.  As $c$ goes from $0$ to $1$, the intersection
cohomology groups isomorphic to the maximal cohomology interpolate between middle perversity
and absolute cohomology, while the intersection cohomology groups isomorphic to the minimal
cohomology interpolate between middle perversity and relative perversity.

Using our Hodge theorem, we can also obtain a signature theorem for the manifolds $(M, g_c)$ through 
a signature theorem for intersection cohomology.  If $p$ and $q$ are dual perversities with $p \geq q$, 
then we can define an intersection form on spaces of the form
$$
\mbox{Im} \, \big(I\!H^{n/2}_{p}(X, B)
\longrightarrow I\!H^{n/2}_{q}(X, B) \big)
$$ 
\noindent
as follows.  Define a (degenerate) pairing on 
$$
I\!H^{k}_{p}(X, B) \otimes I\!H^{n-k}_{q}(X, B)
$$
\noindent via the map induced by inclusion
\begin{equation}
I\!H^{k}_{p}(X, B) \longrightarrow I\!H^{k}_{q}(X, B)
\label{eq:perversityinclusion}
\end{equation}
\noindent
and the nondegenerate pairing 
\begin{equation}
I\!H^{k}_{p}(X, B) \otimes I\!H^{n-k}_{q}(X, B) \longrightarrow \mathbb{R}.
\label{eq:pairing}
\end{equation}
\noindent
Call its signature $\sigma_p(X)$.

The $\tau$-invariant defined by Dai in \cite{Dai} is given by a sum $\tau = \sum_{i=2}^{\infty} \tau_i$, 
where $\tau_i$ is the
signature of a form defined on the $E_i$ term of the Leray-Serre spectral sequence for the boundary
fibration of $M$.  Our signature theorem for metrics interpolating between finite cylindrical and conical is:

\begin{theorem}  If $p = \underline{\frakm} + k$, then the signature of the intersection form on these 
spaces is
given by:
$$
\sigma_p(X)= \mbox{sgn}\,
\mbox{Im}\,\big( H^*(M, \del M) \longrightarrow H^*(M))\big) + \sum_{i=2+2k}^{\infty} \tau_i.
$$
\label{th:signature}
\end{theorem}
\noindent

Thus as the metric becomes less and less cylindrical and more and more conical, 
the signature theorem picks up more and more of the $\tau_i$ terms, until when the 
metric is close to conical, the signature includes all of $\tau$.

The author would like to thank Tam\'{a}s Hausel, Xianzhe Dai, and Rafe Mazzeo for 
useful conversations relating to this paper.

\section{Background}
Before we begin the proof of Theorem \ref{th:hodge}, we will review briefly some definitions and 
theorems we will use in our proofs.

\subsection{Intersection Cohomology}
We will prove $H^*_{max}(M, g_c)$ is isomorphic
to an intersection cohomology group via a sheaf argument using the following
theorem of Goresky and MacPherson:
\begin{proposition}[\cite{GM2}] Let $X$ be a stratified space
and let $(\mathcal{L}^*, d)$ be a complex of fine sheaves on $X$ with
cohomology $H^*(X,\calL)$. Suppose that if $\calU$ is a neighbourhood
in the principal (smooth) stratum of $X$, then $H^*(\calU,\calL) =
H^*(\calU,\CC)$, while if $q$ lies in a stratum of codimension
$\ell$, and $\calU = \calV \times C(L)$ where $\calV$ is a neighborhood of $q$ in the 
codimension $l$ stratum, then
\begin{equation}
H^k(\calU,\calL) \cong I\!H_\frakp^k(\calU) =
\left\{ \begin{array}{ll}
I\!H_\frakp^k(L) & k \leq \ell-2-p(\ell) \\
0         & k \geq \ell-1-p(\ell)
\end{array} \right.
\label{eq:localcalc}
\end{equation}
Then there is a natural isomorphism between the hypercohomology
$\mathbb{H}^{\,*}(X,\mathcal{L}^*)$ associated to this complex
of sheaves and $I\!H^*_\frakp(X)$, the intersection cohomology of
perversity $\frakp$.
\label{pr:shchar}
\end{proposition}

For a pseudomanifold $X$, such as ours, which has only one singular 
stratum, $B$, the local calculation for intersection cohomology in a neighborhood of a point $b$
on the singular stratum
$\calU \cong \Delta \times C(F)$ for $b \in \Delta$ a disk in $B$ is:
\begin{equation}
H^i(\calU,\calL) \cong I\!H_\frakp^i(\calU) =
\left\{ \begin{array}{ll}
H^i(F) & k \leq f-1-p \\
0         & k \geq f-p
\end{array} \right.
\label{eq:localcalc2}
\end{equation}

Note that for larger values of $p$, the truncation in $H^i(F)$ occurs at a lower degree.  Thus if $p>q$,
then $I\!H^i_p(U)$ vanishes for more degrees than $I\!H^i_q(U)$.  

We generalize this definition slightly as follows and use this definition in both theorems:
\begin{equation}
I\!H^*_{j}(X,B) = \left\{ \begin{array}{ll}
H^*(X-B) & j \leq -1, \\
I\!H_\frakp^*(X) & 0 \leq j \leq f-1, \\
H^*(X, B) & j \geq f,
\end{array} \right. .
\label{eq:extih} 
\end{equation}

There are two particular perversities on a pseudomanifold, $X$, called the {\em middle perversities}.
In the case that the strata of $X$ are all even dimensional, $X$ is a particular kind of pseudomanifold 
called a Witt space.  In this case, the two middle perversity intersection cohomologies are identical, and 
this single middle perversity intersection cohomology satisfies Poincar\'{e} duality.  If the strata are not 
all even dimensional, the two middle perversities are not generally the same.  In this case, they are 
Poincar\'{e} duals of each other. Perversities which lie at equal distances from these two middle
perversities, $\overline{\frakm} - k$ and $\underline{\frakm} + k$ are also dual perversities, that is, 
their corresponding intersection cohomologies are Poincar\'{e} duals.

The local calculations near a point on the compactifying layer $B$ for the two middle perversities are:
\begin{equation}
H^k(\calU,\calL) \cong I\!H_{\underline{\frakm}}^k(\calU) =
\left\{ \begin{array}{ll}
H^k(F) & k \leq \frac{f}{2} - 1 \mbox{ if $f$ is even and }  \frac{f-1}{2}  \mbox{ if $f$ is odd} \\
0         & k \geq \frac{f}{2} \mbox{ if $f$ is even and }  \frac{f+1}{2}  \mbox{ if $f$ is odd} 
\end{array} \right.
\end{equation}
for the lower middle perversity and
\begin{equation}
H^k(\calU,\calL) \cong I\!H_{\overline{\frakm}}^k(\calU) =
\left\{ \begin{array}{ll}
H^k(F) & k \leq \frac{f}{2}  \mbox{ if $f$ is even and }  \frac{f-1}{2}  \mbox{ if $f$ is odd} \\
0         & k \geq \frac{f}{2} + 1 \mbox{ if $f$ is even and }  \frac{f+1}{2}  \mbox{ if $f$ is odd} 
\end{array} \right.
\end{equation}
for the upper middle perversity.  Note that if we allow for the possibility that the fibre is a point, that is, if 
our manifold $(M, g_c)$ is a simple manifold with boundary, then these local calculations correspond to 
$H^*_0(M)$ and $H^*(M)$, respectively.  Thus our extended definition of intersection cohomology 
allows us to state the 
results for the case where $F$ is trivial in the same manner as the cases where $F$ is nontrivial.  

\subsection{Geometry of Fibre Bundles}
It is useful to think of the Leray-Serre spectral sequence of a fibration in terms of bidegree.
So recall from geometry of fibre bundles that we can define the bidegree of a form on the 
total space of a fibration as follows.
Let $\phi:Y \to B$ be a fibration with fibre $F$, and suppose that it is endowed
with a metric $G$ of the form $\phi^*(h) + k$, where $h$ is a metric on $B$. We
assume furthermore that $\phi: (Y,G) \to (B,h)$ is a Riemannian submersion.
The tangent bundle $TY$ splits into a vertical and horizontal subbundle, 
$T^V Y \oplus T^H Y$, where $T^V Y =\mbox{ker\,}(d\phi)$ and $T^H Y$ is its
orthogonal complement (and also the subbundle annihilated by $k$). This induces
a splitting of the form bundles on $Y$, and thus every differential form has a 
(horizontal,vertical) bidegree, i.e.
\begin{equation}
\Omega^{p,q}(Y) = \Omega^p(B)\, \widehat{\otimes} \, \Omega^q(Y,T^VY).
\label{eq:decomp}
\end{equation}

\subsection{Hilbert Complexes}
We will use the concept of Hilbert complexes and results about them from \cite{BL} in our Hodge 
result, so we summarize them briefly here.
Consider a complex of the form
\begin{equation}
0 \to L_0 \stackrel{D_0}{\longrightarrow} L_1 
\stackrel{D_1}{\longrightarrow}L_2 \ldots \stackrel{D_{n-1}}{\longrightarrow} L_n \to 0,
\label{eq:hc}
\end{equation}
where each $L_i$ is a separable Hilbert space, $D_i:L_i \to L_{i+1}$ is
a closed operator with dense domain $\calD(D_i)$ such that $\ran(D_i) \subseteq 
\calD(D_{i+1})$ and $D_{i+1}\circ D_i = 0$ for all $i$. Under these
conditions, (\ref{eq:hc}) is called a Hilbert complex, and is denoted by
$(L_*,D_*)$.  

Many familiar constructions in Hodge-de Rham theory carry over
immediately to this setting.  In particular:
\begin{itemize}
\item[i)] There is a dual Hilbert complex 
\[
0 \to L_0 \stackrel{D_0^*}{\longleftarrow} L_1 
\stackrel{D_1^*}{\longleftarrow}L_2 \ldots \stackrel{D_{n-1}^*}{\longleftarrow} L_n \to 0
\]
defined using the Hilbert space adjoints of the differentials, $D_i^*: L_{i+1} \to L_i$;
\item[ii)] The Laplacian $\Delta_i = D_{i}^* D_i + D_{i-1}D_{i-1}^*$
is a self-adjoint operator on $L_i$ with domain
\[
\calD(\Delta_i) = \{u \in \calD(D_i) \cap \calD(D_{i-1}^*): D_i u \in \calD(D_i^*),
D_{i-1}^* u \in \calD(D_{i-1})\} 
\]
and nullspace
\[
\ker \Delta_i := \calH^i(L_*,D_*) = \ker D_i \cap \ker D_{i-1}^*;
\]
\item[iii)] (\cite{BL}, Lemma 2.1) There is a weak Kodaira decomposition
\[
L_i = \calH^i \oplus \overline{\ran D_{i-1}} \oplus \overline{\ran{D_i^*}};
\]
\item[iv)] (\cite{BL}, Corollary 2.5) The cohomology of $(L_*,D_*)$ is defined by
\[
H^i(L_*,D_*) = \ker D_i/\ran D_{i-1};
\]
if this space is finite dimensional, then $\ran D_{i-1}$ is necessarily
closed and $H^i(L_*,D_*) = \calH^i(L_*,D_*)$. 
\item[v)] (\cite{BL}, Corollary 2.6)  The cohomology of the complex $(L_*,D_*)$
is finite dimensional if an only if the cohomology of the dual complex $(L^*_*,D^*_*)$ is.
\end{itemize} 

We will also simplify our calculations using the K\"unneth theorem:
\begin{proposition}[\cite{BL} Corollary 2.15]
Let $(L',D')$ and $(L'' D'') $ be two Hilbert complexes.  Form the 
completed tensor product Hilbert complex $(L,D)$:
\[
L_j = \bigoplus_{i+\ell=j} L'_i \, \hat \otimes \, L''_\ell,
\]
\[
D_j = \bigoplus_{i+\ell=j} (D'_i \otimes \mbox{id}_{L''_\ell} + (-1)^i 
\mbox{id}_{L'_i} \otimes D"_\ell).
\]
\noindent
Suppose that $D''$ has closed range in all degrees. Then
\[
H^j(L,D)= \bigoplus_{i+\ell=j} H^i(L',D') \otimes H^\ell(L''D'').
\]
\label{pr:kunneth}
\end{proposition}

The last result we will use says that we may compute these 
cohomology groups using a `core subcomplex' 
\[
\calD^\infty_{{\max}} \Omega^*(M,g) \subset L^2 \Omega^*(M,g)
\]
consisting of all elements $\omega$ which are in the domain of $\Delta_{M,m}^\ell$ 
for every $\ell \geq 0$.   
\begin{proposition}[\cite{BL} Theorem 2.12]
The cohomology $H^*_{{\max}}(M,g)$ is equal to the cohomology of the
complex $(\calD^\infty \Omega^*_{{\max}}(M,g), d_{{\max}})$.
\label{pr:core}
\end{proposition}

This implies that  $H^*(M)$ can be computed using the complex
of forms which are smooth on the interior of $M$.

\section{Proof of Theorem \ref{th:hodge}}

The two complexes $(D(d_{max}), d_{max})$ and $(D(d_{min}), d_{min})$ form Hilbert complexes as
defined in \cite{BL}, so we can use the theory developed in that paper to prove Theorem \ref{th:hodge}.
First we note that it suffices to prove the theorem for $(D(d_{max}), d_{max})$.  This is because
by Corollary 2.6 of \cite{BL}, the cohomology of this complex is finite dimensional if and only 
if the cohomology of the dual complex $(D(\delta_{min}), \delta_{min})$ is.  The Hodge star operator 
gives a isomorphism of complexes between this dual complex and the complex 
$(D(d_{min}), d_{min})$, so the maximal complex has finite dimensional cohomology if and only if 
the minimal complex does.  These two are Poincar\'{e} dual spaces, so if $H^*_{max} \cong
I\!H^*_{\mathfrak{p}}(X, B)$, then we must also have $H^*_{min} \cong I\!H^*_{\mathfrak{q}}(X, B)$, where
$\mathfrak{p}$ and $\mathfrak{q}$ are dual perversities.  

Further, if these cohomologies are finite dimensional, for instance, if they are isomorphic
to intersection cohomologies, then by Corollary 2.5 in \cite{BL}, we immediately get the isomorphisms
$$
H^*_{min}(M, g_c) \cong  \mathcal{H}^*_{m,M}(M, g_c)
$$
and
$$
H^*_{max}(M, g_c) \cong  \mathcal{H}^*_{M, m}(M, g_c).
$$

We can create a complex of fine sheaves of $L^2$ forms on $X$ in the same manner as in
the proofs of Proposition 2 in \cite{HHM} and Proposition 1.18 in \cite{Z2}.  By Proposition \ref{pr:shchar}, 
we need only then prove a local Poincar\'{e} lemma to establish
the isomorphism between $H^*_{max}(M, g_c)$ and intersection cohomology of some perversity.  
The local Poincar\'{e} lemma for points on the stratum $M \subset X$ is identical to the standard
Poincar\'{e} lemma for compact manifolds.  Near a point on the singular stratum, we need to establish a 
Poincar\'{e} lemma for neighborhoods of the form $U = \Delta \times C(F)$, where $\Delta$ is a disk in
the singular stratum.  By the Kunneth Theorem, Corollary 2.15,  from \cite{BL}, this calculation
reduces:
$$
H^*_{max}(U, g_c) \cong H^*_{max}(\Delta) \otimes H^*_{max}(C(F), g_c) \cong H^*_{max}(C(F), g_c),
$$
\noindent
so the necessary Poincar\'{e} lemma is one for the cone on $F$ with respect to the metric
$g_c = dr^2 + r^{2c}ds_F^2$.

\subsection{Preliminary propositions}

The proof of the necessary Poincar\'{e} lemma is based on the techniques in \cite{Ch}.  
Before proceeding to the proof, we'll lay out a few facts we will use.  
First, by Theorem 2.12 of \cite{BL}, we know that the natural inclusion of smooth $L^2$ forms into 
the space of $L^2$ forms induces an isomorphism on cohomology for both the maximal and 
minimal complexes, so we may always assume a max cohomology class has a smooth 
representative.  This is true for any smooth manifold.  Next we prove a few propositions specific to 
the metrics $g_c$ on cones.  

\begin{proposition}
Let $\phi$ be an $i$-form on $F$.  Then the pullback of $\phi$ to $C(F)$ is in 
$L^2\Omega^i(C(F),g_c)$ if and only if $i < \frac{f}{2} + \frac{1}{2c}$, and in this case,
the pullback map is bounded.
\label{prop:pullbacks}
\end{proposition}
\begin{proof}
If $\phi$ is an $i$-form on $F$, then 
$$
\int_{C(F)} ||\phi||^2_c dvol_c = ||\phi||^2_{L^2(F)} \int_0^1 r^{c(f-2i)} dr <\infty
$$
if and only if $c(f-2i)> -1$, which is if and only if $i < \frac{f}{2} + \frac{1}{2c}$.
Since the integral on the right is independent of $\phi$, the pullback map is bounded.
\end{proof}

\begin{proposition}
There exists a $K>0$ such that for all 
$\alpha = \phi + dr \wedge \omega \in L^2 \Omega^i(C(F), g_c)$ and for any null set $S \subset (1/2, 1)$
there is an $a \in (1/2, 1) - S$ such that 
$$
||\phi(a)||^2_{L^2(F)} \leq K ||\phi||^2_{L^2(C(F), g_c)} \leq K ||\alpha||^2_{L^2(C(F), g_c)}.
$$
\label{prop:bound}
\end{proposition}
\begin{proof}
Suppose not.  Then for any $K>0$, there is some family of forms on $F$, 
$\phi(r)  \in L^2 \Omega^i(C(F), g_c)$ such that 
$$
||\phi||^2_{L^2(C(F), g_c)} \geq \int_{1/2}^1 \int_F ||\phi||^2_F r^{-2ic} dvol_F r^{cf} dr
$$
$$
= \int_{1/2}^1 ||\phi(r)||^2_{L^2(F)} r^{c(f-2i)} dr > K ||\phi||^2_{L^2(C(F), g_c)}
 \int_{(1/2, 1) - S} r^{c(f-2i)} dr
 $$
 $$
 =K ||\phi||^2_{L^2(C(F), g_c)} \int_{1/2}^1 r^{c(f-2i)} dr.
 $$
So choose $K > (\int_{1/2}^1 r^{c(f-2i)} dr)^{-1}$ and we have a contradiction.
\end{proof}

\begin{proposition} If $i < \frac{f}{2} + \frac{1}{2c} +1$ and $\alpha = \phi + dr \wedge \omega \in L^2 
\Omega^i(C(F), g_c)$, then for any $a\in (1/2, 1)$, 
$$
K_a(\alpha) = \int_a^r \omega(s) ds \in L^2 \Omega^{i-1}(C(F), g_c).$$
\label{prop:K}
\end{proposition}
\begin{proof}
This proposition follows essentially from Schwartz's inequality.  
\begin{equation}
||K_a||^2_{L^2(C(F),g_c)} = \int_0^1 \int_F ||\int_a^r \omega(s) \, ds||^2_{F} r^{c(f-2i+2)} \, dvol_F \,dr
\label{eq:Ka}
\end{equation}
since $\omega$ is a family of $i-1$ forms on $F$.  Consider just the inside of this:
$$
||\int_a^r \omega(s) \, ds||^2_{F} \leq ( \int_a^r ||\omega||_F)^2
$$
and for any $j$, by the Schwartz inequality,
$$
=(\int_a^r s^{-2j} || s^j \omega(s)||_F \, ds)^2 
\leq \int_a^r s^{-2j} \, ds \int_a^r ||s^j \omega(s)||^2_F \, ds
$$
$$
= \left\{ 
\begin{array}{ll}
\frac{r^{1-2j} - a^{1-2j}}{1-2j} \int_a^r||s^j \omega(s)||^2_F \, ds & j \neq 1/2 \\
(\ln(r) - \ln(a)) \int_a^r||s^j \omega(s)||^2_F \, ds & j = 1/2
\end{array}
\right. .
$$
Choose $2j = c(f-2i + 2)$ and replace this expression in equation \ref{eq:Ka} to get 
$$
||K_a||^2_{L^2(C(F),g_c)} \leq 
\left\{ 
\begin{array}{ll}
\int_0^1 \frac{r - a(r/a)^{c(f-2i + 2)}}{1-c(f-2i+2)} \int_a^r|| \omega(s)||^2_F s^{c(f-2i+2)}\, ds \, dvol_F \, dr
& i \neq  \frac{f}{2} + 1 - \frac{1}{2c}\\
\int_0^1 r(\ln(r) - \ln(a)) \int_a^r||\omega(s)||^2_F s^{c(f-2i+2)} \, ds \, dvol_F \, dr&  i =  \frac{f}{2} + 1 - 
\frac{1}{2c}
\end{array}
\right. 
$$
$$
\leq
\left\{ 
\begin{array}{ll}
\int_0^1 \frac{r - a(r/a)^{c(f-2i + 2)}}{1-c(f-2i+2)} \, dr || \omega||^2_{L^2(C(F), g_c)}
& i \neq  \frac{f}{2} + 1 - \frac{1}{2c}\\
\int_0^1 r(\ln(r) - \ln(a)) \, dr || \omega||^2_{L^2(C(F), g_c)}&  i =  \frac{f}{2} + 1 - \frac{1}{2c}
\end{array}
\right. .
$$
Since $a \in (1/2, 1)$, the first integral is uniformly bounded in $a$ for $i \leq \frac{f}{2} + \frac{1}{2c}+1$
and the second integral is also bounded uniformly in $a$.  Thus $K_a$ is a bounded operator.
\end{proof}

\begin{proposition} Let $\rho >0$, and endow $(\rho,1) \times F$ with the metric $g_c$ restricted from
$C(F)$.    Let $\alpha = \phi + dr \wedge \omega  \in L^2\Omega^i(C(F), g_c)$
If $i \geq  \frac{f}{2} + \frac{1}{2c}$, then there exists a sequence $\epsilon_s \rightarrow 0$
such that 
$$\lim_{\epsilon_s \rightarrow 0} ||\phi(\epsilon_s)||^2_{L^2((\rho,1) \times F, g_c)} = 0.
$$
\label{prop:epsilons}
\end{proposition}
\begin{proof}
Since $\alpha \in L^2\Omega^i(C(F), g_c)$, so is $\phi$, so we know that
$$
\int_0^1\int_F ||\phi||^2_F \,dvol_F r^{c(f-2i)} \, dr < \infty.
$$
That is, 
$$
\int_F ||\phi(r)||^2_F \,dvol_F r^{c(f-2i)} \in L^1[(0,1)].  
$$
Thus by Lemma 1.2 in \cite{Ch}, there is a sequence $\epsilon_s \rightarrow 0$ for which
$$
|\int_F ||\phi(\epsilon_s)||^2_F \,dvol_F \epsilon_s^{c(f-2i)} | < \frac{C}{\epsilon_s | \ln(\epsilon_s)|}
$$
for some constant $C>0$.  So we have 
$$
|\int_F ||\phi(\epsilon_s)||^2_F \,dvol_F| < \frac{C \epsilon_s^{c(2i-f)-1}}{|\ln(\epsilon_s)|}.
$$
Since $i \geq  \frac{f}{2} + \frac{1}{2c}$, the right hand side tends to zero as $\epsilon_s \rightarrow 0$.
Thus 
$$
||\phi(\epsilon_s)||^2_{L^2((\rho, 1)\times F, g_c)} = \int_{\rho}^1 \int_F
||\phi(\epsilon_s)||^2_F r^{c(f-2i)} \, dvol_F \, dr 
$$
$$
=||\phi(\epsilon_s)||^2_{L^2(F)} \int_{\rho}^1 r^{c(f-2i)} \, dr \longrightarrow 0
$$
also as $\epsilon_s \rightarrow 0$.

\end{proof}

\subsection{Poincar\'{e} lemma}

The Poincar\'{e} lemma we need to prove for Theorem \ref{th:hodge}, is the following :
\begin{lemma}
Let $(F,g)$ be an $f$-dimensional compact manifold and consider the family of metrics on the cone
over $F$, $C(F)$, given by $ds_c^2 = dr^2 + r^{2c}ds_F^2$ for $0 < c < 1$.  Then
$$
H^i_{max}(C(F), g_c) \cong 
\left\{
\begin{array}{ll}
H^i(N) & i < \frac{f}{2} + \frac{1}{2c} \\
0          & i \geq \frac{f}{2} + \frac{1}{2c}
\end{array} \right.
.
$$
\label{lem:poincare}
\end{lemma}

\begin{proof}
Consider first the case where $ i < \frac{f}{2} + \frac{1}{2c}$.  We want to create a bijective bounded 
linear map 
$$R: H^i_{max}(C(F), g_c) \longrightarrow H^i(F).$$
\noindent
 If $\alpha = \phi + dr \wedge \omega \in L^2 \Omega^i(C(F), g_c)$, then for all but a null set of points in 
$(1/2,1)$, we get
$\phi \in L^2(F)$.  So let $[\alpha] \in H^i_{max}(C(F), g_c)$ and for any such value, $a$, define 
$R([\alpha]) = [\phi(a)]$.  We need to check that this map makes sense, that it is independent of our 
choice of cohomology representative and of a, that it is bijective, and that it is bounded with respect
to the natural norm on cohomology:
$$
|| [\alpha] ||_{L^2H} = \mbox{min}_{\gamma \in [\alpha]} ||\gamma||_{L^2}.
$$

First, the map clearly makes sense, since if $\alpha$ is closed, then 
$d \alpha = d_F \phi + dr \wedge (\phi' - d_N \omega)$, thus $d_F \phi = 0$ for all a.  To show
that the map is well defined, suppose that $\alpha= \phi + dr \wedge \omega$ is a 
smooth representative of $[\alpha]$ and 
$\beta = \psi + dr \wedge \nu$ is any other representative.  Then 
$\phi(a) \in L^2\Omega^i(F)$ for any $a$, so choose any $a$ such that $\psi(a) \in L^2\Omega^i(F)$.
Then there is some form $\eta = \rho + dr \wedge \theta \in L^2\Omega^{i-1}(C(F), g_c)$ such 
that $\beta = \alpha + d_{max} \eta$, so $\psi = \phi + d_F \rho$ for all $a$.  Since $L^2$ cohomology
on a compact manifold is the same as absolute cohomology, this means that even if $\rho(a)$ is not
in $L^2(F)$, there must be some $L^2$ form $\tilde {\rho}_a$ which is and for which 
$\psi(a) = \phi(a) + d_F \tilde{\rho}_a$.  Thus $[\psi(a)] = [\phi(a)]$, so the map $R$ is independent
of the choice of cohomology representatives.  To show that it is independent of the choice of $a$,
it suffices therefore to consider smooth representatives.    For $\alpha$ a representative as before, 
since $\phi(a) \in L^2\Omega^i(F)$ for any $a$, we need to show that $[\phi(a)] = [\phi(b)]$ for
any $a,b \in (1/2, 1)$.  Since $\alpha$ is smooth, 
$
\int_a^b \omega(s) \, ds \in L^2\Omega^{i-1}(F)
$
and
$$
d\int_a^b \omega(s) \, ds  = \int_a^b d_F \omega(s) \, ds = \int_a^b \phi'(s) \, ds
$$
$$
= \phi(b) - \phi(a).
$$
\noindent
Thus $[\phi(a)] = [\phi(b)]$, so the map $R$ is well-defined.

To show the map is linear, let $a$ be a value in $(1/2,1)$ such that $\alpha(a)$, $\beta(a)$ and 
$(\alpha+\beta)(a)$ are all in $L^2\Omega^i(F)$.  Then $R([\alpha]) + R([\beta]) = [\alpha(a)]
+[\beta(a)] = [(\alpha + \beta)(a)] = R([\alpha+\beta])$.  
The map is bounded by proposition \ref{prop:bound} and surjective by proposition \ref{prop:pullbacks}.  
So we have only left to show that 
it is injective.  So suppose that $R([\alpha])=[0]$ and let $\alpha= \phi + dr \wedge \omega$ be a smooth 
representative.
Then $\phi(a) = d_F \eta$ for some $\eta \in L^2\Omega^{i-1}(F)$.  By proposition \ref{prop:pullbacks}, 
we can consider
$\eta$ as a form in $L^2\Omega^{i-1}(C(F), g_c)$, and by proposition \ref{prop:K}, 
$
\int_a^r \omega(s) \, ds \in L^2\Omega^{i-1}(C(F), g_c)
$.
Then 
$$
d(\eta + \int_a^r \omega(s) \, ds) = d_F \eta + dr \wedge \omega(r) + \int_a^r d_F\omega(s)\, ds 
$$
$$
= \phi(a) + dr \wedge \omega(r) + \int_a^r \phi'(s)\, ds = \alpha.
$$
Thus $\alpha$ is exact and $[\alpha] = [0]$.

Now consider the case where $i \geq \frac{f}{2} + \frac{1}{2c}$.  Since $c<1$, this implies that 
$c(f-2i+2) < 1$.  We want to show that any class $[\alpha] \in H^i_{max}(C(F),g_c)$ is trivial. 
First we know again that any such class can be represented by a smooth form.  
Let $\alpha = \phi + dr \wedge \omega \in L^2 \Omega^i(C(F), g_c)$ be a smooth representative
of $[\alpha]$.  We need
to show that $\alpha = d_{max} \eta$ for some $\eta \in L^2 \Omega^{i-1}(C(F), g_c)$.   
Consider $K_0 (\alpha) = \int_0^r \omega(s) \, ds$.  Then 
$$
||K_0 \alpha||^2_{L^2(C(F),g_c)} = \int_0^1 \int_F || \int_0^r \omega(s)\, ds||_F^2 r^{c(f-21 + 2)} \, dvol_F \, 
dr.
$$
By the same steps as in the proof of proposition \ref{prop:K}, we get that this is 
$$
\leq \int_0^1 \int_F r^{c(f-2i+2)} [ \int_0^r s^{-2j} \int_0^r ||s^j \omega(s)||_F^2 \,ds]\, dvol_F \, dr
$$
for all $j$.  For $2j < 1$, this is 
$$
= \int_0^1 \int_F r^{c(f-21+2)} r^{-2j+1} [ \int_0^r ||s^j \omega(s)||_F^2 \, ds] \, dvol_F \, dr.
$$
Let $2j = c(f-2i + 2)$.  Then this becomes
$$
\int_0^1 \int_F r \int_0^r ||\omega(s)||^2_F s^{c(f-2i+2)} \, ds \, dvol_F \, dr
\leq 
\int_0^1r  ||\omega(s)||^2_{L^2(C(F), g_c) }\, dr \leq \frac{1}{2} ||\alpha||^2_{L^2(C(F), g_c) }.
$$
Thus $K_0: L^2 \Omega^i(C(F),g_c) \longrightarrow L^2 \Omega^{i-1}(C(F),g_c)$ is a bounded map.
We want to show that if $\alpha$ is closed, then $d_{max}(K_0(\alpha)) = \alpha$.  
This means that we need to show for any $\beta \in C^{\infty}_0 \Omega^{i-1}(C(F))$, we have
\begin{equation}
\langle K_0(\alpha), \delta \beta \rangle_{L^2(C(F), g_c)} = \langle \alpha, \beta \rangle _{L^2(C(F), g_c)}.
\label{eq:K0dmax}
\end{equation}
So let $\beta \in C^{\infty}_0 \Omega^{i-1}(C(F))$.  Then for some $\rho > 0$, $\beta$ is supported in
$(\rho, 1) \times F$.  So equation \ref{eq:K0dmax} becomes
\begin{equation}
\langle K_0(\alpha), \delta \beta \rangle_{L^2((\rho,1) \times F, g_c)} = \langle \alpha, \beta \rangle 
_{L^2((\rho, 1) \times F, g_c)}.
\label{eq:dmax2}
\end{equation}
We'd like to use the fundamental theorem of calculus as in the proof for $i< \frac{f}{2} + \frac{1}{2c}$, 
but since we don't know that $\alpha$ is continuous at 0, we can't do this directly.  So to get around 
this we use the following proposition:
\begin{proposition} Define 
$$K_{\epsilon} (\alpha) = \int_{\epsilon}^r \omega(s)\, ds.$$
Then on $(\rho , 1) \times F$ with the restricted $g_c$ metric, 
$$
K_{\epsilon}(\alpha) \stackrel{L^2((\delta,1) \times F, g_c)}{\longrightarrow} K_0(\alpha).
$$
\label{prop:Kepsilon}
\end{proposition}
\begin{proof}
We have
$$
||K_{\epsilon}(\alpha) - K_0(\alpha)||^2_{L^2((\delta,1) \times F, g_c)}
= \int_{\delta}^1 \int_F ||\int_0^{\epsilon} \omega(s)\, ds||_F r^{c(f-2i + 2)} \, dvol_F \, dr.
$$
As in the proof of Proposition \ref{prop:K}, this is
$$
 \leq \epsilon^{-c(f-2i + 2) +1} \int_{\delta}^1  r^{c(f-2i + 2)} \, dr ||\int_0^{\epsilon} \omega(s) \, 
ds||^2_{L^2(C(F), g_c)}.
$$
Since $c(f-2i + 2) < 1$, the first integral in this product converges, and the whole expression
tends to 0 as $\epsilon \rightarrow 0$.  
\end{proof}

So now we can consider for a closed smooth $\alpha$, 
$$
dK_{\epsilon} (\alpha) = d \int_{\epsilon}^r \omega(s) \, ds = dr \wedge \omega(r) + 
\int_{\epsilon}^{r} d_F \omega(s) \, ds
$$
$$
= dr \wedge \omega(r) + \int_{\epsilon}^r \phi'(s) \, ds
= dr \wedge\omega(r) + \phi(r) - \phi(\epsilon) = \alpha - \phi(\epsilon).
$$
By proposition \ref{prop:epsilons}, there is some sequence $\epsilon_s \rightarrow 0$ such that
$$\lim_{\epsilon_s \rightarrow 0} ||\phi(\epsilon_s)||^2_{L^2((\rho,1) \times F, g_c)} = 0.
$$
So our equation \ref{eq:dmax2} can be proved as follows:
$$
\langle K_0(\alpha), \delta \beta \rangle_{L^2((\rho,1) \times F, g_c)} = 
\lim_{\epsilon \rightarrow 0}
 \langle K_{\epsilon}(\alpha), \delta \beta \rangle_{L^2((\rho,1) \times F, g_c)} 
 = \lim_{\epsilon \rightarrow 0}
 \langle dK_{\epsilon}(\alpha), \beta \rangle_{L^2((\rho,1) \times f, g_c)} 
$$
$$
= \langle \alpha, \beta \rangle_{L^2((\rho,1) \times F, g_c)} -
 \lim_{\epsilon \rightarrow 0} \langle \phi(\epsilon),  \beta \rangle_{L^2((\rho,1) \times F, g_c)} .
$$
Thus the limit 
$$
 \lim_{\epsilon \rightarrow 0} \langle \phi(\epsilon),  \beta \rangle_{L^2((\rho,1) \times F, g_c)} 
 $$
 must exist.  This means we can calculate it from any subsequence, such as the subsequence
 $\epsilon_s$ given in proposition \ref{prop:epsilons}.  So we get
 $$
  \lim_{\epsilon \rightarrow 0} \langle \phi(\epsilon),  \beta \rangle_{L^2((\rho,1) \times F, g_c)} 
  =
   \lim_{\epsilon_s \rightarrow 0} \langle \phi(\epsilon_s),  \beta \rangle_{L^2((\rho,1) \times F, g_c)} 
   $$
   $$
  \leq  \lim_{\epsilon_s \rightarrow 0} ||\phi(\epsilon_s)||_{L^2((\rho,1) \times F, g_c)} 
  ||\beta||_{L^2((\rho,1) \times F, g_c)} = 0.
  $$
  So 
$$
\langle K_0(\alpha), \delta \beta \rangle_{L^2((\rho,1) \times F, g_c)} 
=  \langle \alpha, \beta \rangle_{L^2((\rho,1) \times F, g_c)} 
$$
and we have that $d_{max} K_0(\alpha) = \alpha$ as required.
\end{proof}


\section{Intersection pairings and Novikov additivity}

In order to prove Theorem \ref{th:signature}, we first need to check that it makes sense to talk about
a signature on Im$(I\!H^{n/2}_p (X) \rightarrow I\!H^{n/2}_q (X) )$, where $p>q$ are dual 
perversities for $X$.  

\begin{lemma} Let $p>q$ be dual perversities for intersection cohomology on the 
compactification $X$ of a manifold with fibration boundary $M$ described in
the introduction.  There is a well defined (degenerate) bilinear pairing 
 $$I\!H^{n/2}_p (X) \otimes I\!H^{n/2}_p (X) \rightarrow \mathbb{R}$$
 given by the map 
 $$I\!H^{n/2}_p (X) \rightarrow I\!H^{n/2}_q (X) )$$
 and the nondegenerate bilinear pairing 
 $$I\!H^{n/2}_p (X) \otimes I\!H^{n/2}_q (X) ) \rightarrow \mathbb{R}.$$
 It descends to a well-defined nondegenerate bilinear pairing $B_X$
 $$ \mbox{Im}(I\!H^{n/2}_p (X) \rightarrow I\!H^{n/2}_q (X) ) \otimes 
 \mbox{Im}(I\!H^{n/2}_p (X) \rightarrow I\!H^{n/2}_q (X) )  \rightarrow \mathbb{R}.
 $$
 \label{lem:intpair}
\end{lemma}
\begin{proof}
By Theorem \ref{th:hodge}, $I\!H^*_p(X)$ and $I\!H^*_q(X)$ can be identified with 
$H^*_{min}(M,g_c)$ and $H^*_{max}(M,g_c)$ for some $0 \leq c \leq 1$.  
To show that the pairing is well defined and descends to one which is well defined on 
$\mbox{Im}(I\!H^{n/2}_p (X) \rightarrow I\!H^{n/2}_q (X))$, 
let $\phi, \psi, \theta \in \mathcal{D}(d_{min})$ be closed
and assume that $[\phi]=[\psi] \in H^*_{max}(M,g_c)$.  Then $\delta_{min} * \theta = 
\pm * d_{min} \theta = 0$ and 
$\phi = \psi + d_{max}\eta$ for some $\eta \in \mathcal{D}(d_{max})$.  So
$$
\int_M \phi \wedge \theta = \int_M (\psi + d_{max} \eta ) \wedge \theta 
= \int_M \psi \wedge \theta + \int_M d_{max}\eta \wedge \theta
$$
$$
= \int_M \psi \wedge \theta \pm \int_M d_{max} \wedge *^2 \theta
= \int_M \psi \wedge \theta \pm \langle d_{max}\eta , *\theta \rangle
$$
$$
= \int_M \psi \wedge \theta \pm \langle \eta , \delta_{min} *\theta \rangle
= \int_M \psi \wedge \theta.
$$
Thus the pairing is well defined and descends to one which is also well defined.
To show it descends to a nondegenerate pairing, we need to show that the subspace
$$
W = \{[\phi] \in H_{min}^{n/2}(M,g_c) | \int_M \phi \wedge \psi = 0 \mbox{ for all }
[\psi] \in H_{min}^{n/2}(M,g_c) \}
$$
is the same as the kernel of $i: H^{n/2}_{min}(M,g_c) \rightarrow H^{n/2}_{max}(M,g_c)$.
Let $[\phi] \in \mbox{Ker}(i)$.  Then $\phi = d_{max} \eta$ for some 
$\eta \in \mathcal{D}(d_{max})$, and for $[\psi]$ in $\mathcal{D}(d_{min})$, we have
$$
\int_M \phi \wedge \psi = \int_M d_{max} \eta \wedge \psi =\pm  \int_M d_{max}\eta  \wedge *^2 \psi
$$
$$
=\pm  \langle d_{max} \eta , * \psi \rangle = \pm \langle \eta , \delta_{min} * \psi \rangle = 0.
$$
Thus $[\phi] \in W$.  Now suppose that $[\phi]$ is not in $ \mbox{Ker}(i)$.  Then
$[\phi] = i([\phi]) \neq [0] \in H^{n/2}_{max}(M,g_c)$.  Thus by Poincar\'{e} duality, 
there is some $[\psi] \in H^{n/2}_{min}(M,g_c)$ with 
$\int_M \phi \wedge \psi \neq 0$.  Thus 
$[\phi]$ is not in $W$.  
So $W$ is exactly Ker$(i)$, and we are done.
\end{proof}

We also need to define a signature on the cone-bundle neighborhood of 
the singular stratum of $X$.  Let $Z$ be a pseudomanifold with boundary,
$Y$, which does not intersect the singular part of $Z$.  
Then define $I\!H^{k}_{p}(Z,Y)$ to be the cohomology of $Z$ with $p$-perversity intersection 
cohomology conditions near the singular set of $Z$ and with relative cohomology conditions 
at the boundary, $Y$.  Define 
$I\!H^{k}_{p}(Z)$ to be the cohomology of $Z$ with $p$-perversity intersection cohomology conditions
near the singular set of $Z$ and with absolute cohomology conditions at the boundary, $Y$.

Define a (degenerate) pairing on
$$ 
I\!H^{k}_{p}(Z,Y) \otimes I\!H^{n-k}_{q}(Z,Y)
$$
\noindent via the map
$$
I\!H^{k}_{p}(Z,Y) \longrightarrow I\!H^{k}_{q}(Z,Y) \longrightarrow I\!H^{k}_{q}(Z) 
$$
\noindent
and the nondegenerate pairing $B_{ZY}$
$$
I\!H^{k}_{p}(Z,Y) \otimes I\!H^{n-k}_{q}(Z) \longrightarrow \mathbb{R}.
$$
\noindent
Call its signature $\hat{\sigma}_p(Z)$.  In the special case that $Z$ is just a manifold with 
boundary, we suppress the $p$, since the signature will be the same for any perversity, 
and will simply be the signature of the intersection pairing on 
$Im(H^{n/2}(Z, \del Z) \rightarrow H^{n/2}(Z))$.
The proof that this pairing is well defined is similar to the proof of Lemma \ref{lem:intpair}.

Next we need a theorem which allows us to study the signature of a neighborhood 
of the boundary separately from the interior of $M$.  We use a version of Novikov Additivity for this.  

\begin{theorem}
If $X$ is a pseudomanifold with a single compact smooth singular stratum and if $Y \subset X$ is a compact codimension 1 submanifold such that $X = Z \bigcup_Y Z^{\prime}$ where
$Z \subset \subset X^{\mbox{reg}}$, then 
$$
\sigma_p(X) = \hat{\sigma}(Z) + \hat{\sigma}_p(Z^{\prime}).
$$
\label{th:Novikov}
\end{theorem}
\begin{proof}
The proof is a modification of the original Novikov additivity theorem in \cite{AS3}.
Assume that $X$ is  $n$ dimensional, and
let $\hat{H}^{n/2}(Z) = Im(H^{n/2}(Z,Y) \rightarrow H^{n/2}(Z)$ and 
$\hat{I\!H}_p^{n/2}(Z^{\prime}) = Im(I\!H_p^{n/2}(Z^{\prime},Y) \rightarrow I\!H_q^{n/2}(Z^{\prime})$

We have exact relative cohomology sequences for $I\!H^*_p(X) $ and $I\!H^*_q(X) $:
$$
\xymatrix{
\cdots \ar[r] & I\!H^{n/2}_p(X,Z)  \ar[r]^{\alpha_p^{\prime}} & I\!H^{n/2}_p(X)  \ar[r]^{\beta_p} 
	& I\!H^{n/2}_p(Z)  \ar[r] & \cdots \\ 
\cdots \ar[r] & I\!H^{n/2}_q(X,Z)  \ar[r]^{\alpha_q^{\prime}} & I\!H^{n/2}_q(X)  \ar[r]^{\beta_q} 
	& I\!H^{n/2}_q(Z)  \ar[r] & \cdots \\ 
\cdots & I\!H^{n/2}_q(Z^{\prime})  \ar[l] & I\!H^{n/2}_q(X)  \ar[l]_{\beta_q^{\prime}} 
	& I\!H^{n/2}_q(X,Z^{\prime})  \ar[l]_{\alpha_q} & \cdots \ar[l]\\ 
\cdots & I\!H^{n/2}_p(Z^{\prime})  \ar[l] & I\!H^{n/2}_p(X)  \ar[l]_{\beta_p^{\prime}} 
	& I\!H^{n/2}_p(X,Z^{\prime})  \ar[l]_{\alpha_p} & \cdots \ar[l]\\ 
}
$$
\noindent
Because $Z$ is smooth, these sequences are the same as:
$$
\xymatrix{
\cdots \ar[r] & I\!H^{n/2}_p(Z^{\prime},Y)  \ar[r]^{\alpha_p^{\prime}} \ar[d]_{i_{Z^{\prime}Y}}
	& I\!H^{n/2}_p(X)  \ar[r]^{\beta_p} \ar[d]_{i_X}
	& H^{n/2}(Z)  \ar[r] \ar[d]_{=}
	& \cdots \\ 
\cdots \ar[r] & I\!H^{n/2}_q(Z^{\prime},Y)  \ar[r]^{\alpha_q^{\prime}} & I\!H^{n/2}_q(X)  \ar[r]^{\beta_q} 
	& H^{n/2}(Z)  \ar[r] & \cdots \\ 
\cdots & I\!H^{n/2}_q(Z^{\prime})  \ar[l] & I\!H^{n/2}_q(X)  \ar[l]_{\beta_q^{\prime}} 
	& H^{n/2}(Z,Y)  \ar[l]_{\alpha_q} & \cdots \ar[l]\\ 
\cdots & I\!H^{n/2}_p(Z^{\prime})  \ar[l] \ar[u]^{i_{Z^{\prime}}}
	& I\!H^{n/2}_p(X)  \ar[l]_{\beta_p^{\prime}} \ar[u]^{i_X}
	& H^{n/2}_p(Z,Y)  \ar[l]_{\alpha_p} \ar[u]^{=}& \cdots \ar[l] \\ 
}
$$
\noindent
The maps $i_{Z^{\prime}Y}$, $i_{Z^{\prime}}$ and $i_X$ are induced from the natural inclusion maps on the complexes of forms, and this entire diagram commutes.  

Define 
$$
A^{\prime} = Im(i_X \circ \alpha_p^{\prime}) \subset Ker(\beta)
$$
$$
A = Im (\alpha_q) = Ker(\beta_q^{\prime} )= Im(i_x \circ \alpha_p).
$$
\begin{proposition}
The spaces $A$ and $A^{\prime}$ are mutual annihilators in 
$Im(I\!H^{n/2}_p(X)  \rightarrow I\!H^{n/2}_q(X))$ under the intersection pairing $B_x$.
\end{proposition}
\begin{proof}
Suppose that $i_X \circ \alpha_p^{\prime}[a] \in A^{\prime}$ and $\alpha_q[b] \in A$.
Then 
$$
B_X(i_X \circ \alpha_p^{\prime}[a] , \alpha_q[b]) = B_{Z^{\prime}Y}([a], bq^{\prime} \circ \alpha_q[b])
=0.
$$
\noindent
So $A$ and $A^{\prime}$ are contained in each other's annihilators.  

Now suppose that $[e] \in Im(I\!H_p^{n/2}(X) \rightarrow I\!H_q^{n/2}(X))$ annihilates $A^{\prime}$
under $B_X$.  Then for every $[c^{\prime}] \in I\!H_p^{n/2}(Z^{\prime},Y)$, 
$$
0 = B_X(i_X \circ \alpha_p^{\prime}[c^{\prime}],[e]) = B_{Z^{\prime}Y}([c^{\prime},\beta_q^{\prime}[e]).
$$
\noindent
Since the pairing $B_{Z^{\prime}Y}$ is nondegenerate, this can only happen if 
$\beta_q^{\prime}[e]=0$.  But then by exactness of the $I\!H_q$ sequence we have 
$[e] \in Im(\alpha_q) = A$.  So $A$ is the annihilator of $A^{\prime}$ under $B_Z$.

Finally, suppose that $i_X[e^{\prime}] \in Im(I\!H_p^{n/2}(X) \rightarrow I\!H_q^{n/2}(X))$ annihilates
$A$ under $B_X$.  Then for all $[c] \in H^{n/2}(Z,Y)$, 
$$
0 = B_X( \alpha_q[c], i_X[e^{\prime}]) = B_{ZY}([c], \beta_q\circ i_X[e^{\prime}])
= B_{ZY}([c], \beta_p[e^{\prime}]).
$$
\noindent
This implies that $\beta_p[e^{\prime}]=0$, so by exactness of the $I\!H_p$ sequence, 
$[e^{\prime}] \in Im (\alpha_p^{\prime})$.  So $i_X[e^{\prime}] \in A^{\prime}$,
and we get that $A^{\prime}$ is also the annihilator of $A$.
\end{proof}

Since the pairing $B_Z$ is nondegenerate, we get as in the original Novikov additivity proof
that 
$$
(A \cap A^{\prime})^* \cong \frac{Im(I\!H_p^{n/2}(X) \rightarrow I\!H_q^{n/2}(X))}{A+A^{\prime}}.
$$
\noindent
Also, 
$$
\frac{A+A^{\prime}}{A \cap A^{\prime}} = \frac{A}{A \cap A^{\prime}} \oplus \frac{A^{\prime}}{A \cap A^{\prime}}
\cong
\frac{Im(\alpha_q)}{Im(\alpha_q) \cap Im(i_X \circ \alpha_p^{\prime})} \oplus
\frac{Im(i_X \circ \alpha_p^{\prime})}{Im(\alpha_q) \cap Im(i_X \circ \alpha_p^{\prime})}.
$$
\noindent
From the exact sequences, we know that
$$
\frac{Im(i_X \circ \alpha_p^{\prime})}{Im(\alpha_q) \cap Im(i_X \circ \alpha_p^{\prime})} \cong
\frac{Im(i_X \circ \alpha_p^{\prime})}{ker(\beta_q^{\prime}) \cap Im(i_X \circ \alpha_p^{\prime})}
\cong Im(\beta_q^{\prime} \circ i_z \circ \alpha_p^{\prime}) \cong 
\hat{I\!H}^{n/2}_p(Z^{\prime}).
$$
\noindent
And also that
$$
\frac{Im(\alpha_q)}{Im(\alpha_q) \cap Im(i_X \circ \alpha_p^{\prime})} 
\cong
\frac{Im(i_X \circ \alpha_p)}{Im(i_X \circ \alpha_p)  \cap i_X(Ker(\beta_q))}
\cong \frac{Im(i_X \circ \alpha_p)}{Im(i_X \circ \alpha_p)  \cap Im(i_X) \cap (Ker(\beta_q)}.
$$
\noindent
Since $Im(i_X \circ \alpha_p) \subset Im(i_X)$, we can eliminate the $Im(i_X)$ term in 
the denominator to get
$$
\cong \frac{Im(i_X \circ \alpha_p)}{Im(i_X \circ \alpha_p) 
\cap Ker(\beta_q)}
\cong 
\frac{Im(\alpha_q)}{Im(\alpha_q) \cap Ker(\beta_q)} \cong \hat{H}^{n/2}(Z).
$$
\noindent
So altogether, we have
$$
\frac{A+A^{\prime}}{A \cap A^{\prime}} \cong \hat{H}^{n/2}(Z) \oplus \hat{I\!H}^{n/2}_p(Z^{\prime}).
$$

By splitting the maps
$$
A \cap A^{\prime} \rightarrow A + A^{\prime} \rightarrow Im(I\!H_p^{n/2}(X) \rightarrow I\!H_q^{n/2}(X))
$$
\noindent
we get that 
$$
Im(I\!H_p^{n/2}(X) \rightarrow I\!H_q^{n/2}(X)) \cong  \hat{H}^{n/2}(Z) \oplus \hat{I\!H}^{n/2}_p(Z^{\prime})
\oplus [(A \cap A^{\prime}) \oplus (A \cap A^{\prime})^*].
$$
As in the proof of the original Novikov additivity, by choosing a good splitting, we can arrange
for the form $B_X$ to be given with respect to this splitting by the block matrix
$$
\left(
\begin{array}{lll}
B_Z & 0 & 0 \\
0 & B_{Z^{\prime}} & 0 \\
0 & 0 & C
\end{array}
\right),
$$
\noindent
where $C$ is the natural form on $[(A \cap A^{\prime}) \oplus (A \cap A^{\prime})^*]$.  The 
signature of such a form is always 0 (p. 580 \cite{AS3}), so we obtain 
$$
\sigma_p(X) = \hat{\sigma}_p(Z) + \hat{\sigma}_p(Z^{\prime})
$$
\noindent
as we wanted.

\end{proof}


\section{Signature of the end}

Now return to our original situation, where $\overline{M}$ is an $n$ dimensional manifold 
with boundary $Y\stackrel{\phi}{\rightarrow} B$
with fibre $F$ and $X$ is the pseudomanifold formed by coning off the fibres of $Y$.  As before,
let $f$ be the dimension of F and $b$ be the dimension of $B$ so that $f+b+1=n$.
Decompose $X$ as $X^{reg} \cup_Y C_{\phi}Y$, where $C_{\phi}Y$ is the pseudomanifold with 
boundary formed by coning off the fibres of $Y$.  Assume that $p = \underline{\frakm} + k$ and
$q = \overline{\frakm} -k$.  
By Theorem \ref{th:Novikov}, to prove Theorem \ref{th:signature}, it suffices to prove the
\begin{lemma}
The signature on $C_{\phi}Y$ is 
$$
\hat{\sigma}_p(C_{\phi}Y^{\prime}) = - \sum_{i= 2 + 2k} ^{\infty}\tau_i.
$$
\label{lem:signend}
\end{lemma}
\noindent
Once we have this, reversing orientation to glue, we get
$$
\sigma_p(X) = \hat{\sigma} (M) + \hat{\sigma}_p(C_{\phi}Y) = 
sgn \, Im(H^{n/2}(M, \del M) \rightarrow H^{n/2}(M)) + \sum_{i= 2 + 2k} ^{\infty}\tau_i
$$
\noindent as required.  

The proof of Lemma \ref{lem:signend} generalizes the techniques of \cite{CD}.  
\begin{proof}
The first step is to identify $I\!H_q^*(C_{\phi}Y)$ and $Im(I\!H_p^*(C_{\phi}Y, Y) \rightarrow I\!H_q^*(C_{\phi}Y))$ in terms of the Leray-Serre spectral sequence for the fibration on $Y$.
We will prove the 
\begin{proposition} In terms of the $E_\infty^{i,j}(Y)$ terms of the Leray-Serre spectral sequence and 
its differentials $d_r^{i,j}: E_r^{i,j}(Y) \rightarrow E_r^{i+r, j-r+1}(Y)$, we get the isomorphisms:
$$
I\!H_q^{n/2}(C_{\phi}Y) \cong \left\{
\begin{array}{ll}
\bigoplus_{\stackrel{i+j=n/2}{j\leq f/2 + k}} (E_{\infty}^{i,j}(Y) \oplus \sum_{s=1}^{\infty} Im(d_s^{i-s, j+s-1})
& b/2 - k >0 \\
\bigoplus_{\stackrel{i+j=n/2}{ j\leq f/2 + k}} E_{\infty}^{i,j}(Y) & b/2 - k \leq 0
\end{array} \right.
$$
$$
Im(I\!H_p^{n/2}(C_{\phi}Y, Y) \rightarrow I\!H_q^{n/2}(C_{\phi}Y)) \cong  \left\{
\begin{array}{ll}
\bigoplus_{\stackrel{i+j=n/2}{ j \geq (f+3)/2 + k}}  \sum_{t \geq j-(f-1)/2 +k} Im(d_t^{i, j-1})
& b/2 - k >0 \\
0 & b/2 - k \leq 0
\end{array} \right.
$$
\end{proposition}
\begin{proof}
Recall that by Theorem \ref{th:hodge}, $I\!H_q^*(C_{\phi}Y) \cong H_{max}^*(\CY, g_c)$
for some metric $g_c$, $0\leq c \leq 1$ 
and $I\!H_p^*(C_{\phi}Y) \cong H_{min}^*(\CY, g_c)$, where $H_{min}^*(\CY, g_c)$ means 
cohomology which is  minimal with respect to the metric $g_c$ near the singular stratum of
$\CY$ and absolute near the boundary $Y$.  Further, as for regular cohomology, 
$I\!H^*_p(C_{\phi}Y,Y)$ is the cohomology of the complex
$(\Omega^{*-1}(Y) \oplus \Omega_{min}^*(C_{\phi}Y, g_c), d_{rel})$ where $\Omega_{min}^*(C_{\phi}Y)$
is the space of smooth forms on $C_{\phi}Y$ which satisfy minimal boundary conditions
for $g_c$ near the singular stratum and absolute boundary conditions
near the boundary $Y$ and where $d_{rel}(\theta, \alpha) = (d_Y \theta - \alpha(1), d_{c} \alpha)$.  
So the short exact sequence
$$
0 \rightarrow \Omega^{*-1}(Y) \stackrel{\gamma}{\rightarrow} \Omega^{*-1}(Y) \oplus \Omega^*_{min}(\CY)
\stackrel{\beta}{\rightarrow} \Omega^*_{min}(\CY) \rightarrow 0
$$
\noindent
induces a long exact sequence on cohomologies.  Similarly, we get a long exact sequence on 
cohomologies corresponding to the perversity $q$:
$$
0 \rightarrow \Omega^{*-1}(Y) \stackrel{\gamma}{\rightarrow} \Omega^{*-1}(Y) \oplus \Omega^*_{max}(\CY)
\stackrel{\beta}{\rightarrow} \Omega^*_{max}(\CY) \rightarrow 0
$$
Using the metric $g_c$, each of these spaces of forms can be filtered by 
(conic) fibre and base bi-degree as in section ?? of \cite{HHM}.
For any any base degree, {i}, we get the following commutative diagram
relating $q$-perversity complexes:
\begin{equation}
\begin{xy}
\xymatrix{
&\vdots &  \vdots  & \vdots & \\
0 \ar[r] &
\Omega^{i,1}(Y) \ar[u]^{d_F} \ar[r]^{\gamma} &
\Omega^{i,1}(Y) \oplus \Omega^{i,2}_{max}(\CY) \ar[u]^{d_{rel}} \ar[r]^{\beta} &
 \Omega^{i,2}_{max}(\CY) \ar[u]^{d_c} \ar[r] & 0\\
0 \ar[r] &
\Omega^{i,0}(Y) \ar[u]^{d_F} \ar[r]^{\gamma} &
\Omega^{i,0}(Y) \oplus \Omega^{i,1}_{max}(\CY) \ar[u]^{d_{rel}} \ar[r]^{\beta}  &
 \Omega^{i,1}_{max}(\CY) \ar[u]^{d_c} \ar[r]  & 0\\
&
0 \ar[u] \ar[r] &
0 \oplus \Omega^{i,0}_{max}(\CY) \ar[u]^{d_{rel}} \ar[r]^{\beta} &
 \Omega^{i,0}_{max}(\CY) \ar[u]^{d_c} \ar[r] & 0
}
\end{xy}
\label{eq:fibreexactsequences}
\end{equation}
The $E_1$ terms of the Leray-Serre spectral sequence for each base degree $i$ fit into the exact sequence:
$$
\cdots \rightarrow E_1^{i, j-1}(\CY,q) \rightarrow
E_1^{i,j-1}(Y) \rightarrow E_1^{i,j}(\CY, Y,q) \rightarrow
E_1^{i, j}(\CY,q) \rightarrow E_1^{i,j}(Y) \rightarrow\cdots .
$$
\noindent
By the Poincar\'{e} lemma for the maximal complex, Lemma \ref{lem:poincare} we get  that for $j \leq \frac{f-1}{2} +k$ if $f$ is odd and $j \leq \frac{f}{2} +k$ if $f$ is even, this is
$$
\rightarrow \Omega^i(B, H^{j-1}(F))
\rightarrow \Omega^i(B,H^{j-1}(F)) \rightarrow E_1^{i,j}(\CY, Y,q) \rightarrow
\Omega^i(B, H^j(F)) \rightarrow  \Omega^i(B,H^{j}(F)) \rightarrow .
$$
\noindent
So for $j \leq \frac{f-1}{2} +k$ if $f$ is odd and $j \leq \frac{f}{2} +k$ if $f$ is even, we find 
$E_1^{i,j}(\CY, Y,q)  =0$.
For $j \geq \frac{f-1}{2} +k +2$ if $f$ is odd and $j \geq \frac{f}{2} +k+2$ if $f$ is even, we have
$$
\cdots \rightarrow 0
\rightarrow \Omega^i(B,H^{j-1}(F)) \rightarrow E_1^{i,j}(\CY, Y,q) \rightarrow
0 \rightarrow \Omega^i(B,H^{j}(F)) \rightarrow\cdots .
$$
\noindent
So for $j \geq \frac{f-1}{2} +k +2$  if $f$ is odd and $j \leq \frac{f}{2} +k+2$ if $f$ is even, we find $E_1^{i,j}(\CY, Y,q) \cong \Omega^i(B,H^{j-1}(F))$.
Thus for the boundary degree $j= \frac{f-1}{2} +k+1$ if $f$ is odd and $j = \frac{f}{2} +k+1$ if $f$ is even, we have (dropping the last term and adding one previous term):
$$
\rightarrow 0 \rightarrow \Omega^i(B, H^{(f-1)/2 + k}(F))
\rightarrow \Omega^i(B,H^{(f-1)/2 +k}(F)) \rightarrow E_1^{i,(f-1)/2+k+1}(\CY, Y,q) \rightarrow
0 \rightarrow ,
$$
\noindent
so $E_1^{i,(f-1)/2+k+1}(\CY, Y,q)=0$ if $f$ is odd and $E_1^{i,f/2+k+1}(\CY, Y,q)=0$ if $f$ is even.  

The differential $d_1$ is horizontal, so it does not
see the truncations.  Thus the $E_2$ levels of the relative and absolute $p$-perversity spectral sequences are for $f$ odd:

Relative $q$ sequence:
\begin{equation}
\begin{array}{llllll}
j=f+1 & &H^0(B, H^{f}(F)) & H^1(B, H^{f}(F)) & H^1(B, H^{f}(F)) 	&\cdots \\
\vdots && \vdots                 & \vdots                  & \vdots                    & \cdots \\
j = \frac{f+3}{2} +k & & H^0(B, H^{(f+1)/2 +k}(F)) & H^1(B, H^{(f+1)/2 +k}(F)) & H^2(B, H^{(f+1)/2 +k}(F)) 	
	&\cdots \\
j = \frac{f+1}{2} +k & & 0&0&0&\cdots \\
\vdots && \vdots & \vdots & \vdots & \cdots \\
j =0 & & 0&0&0&\cdots \\
E_2^{i,j}(\CY, Y, q) && i=0&i=1&i=2& \cdots
\end{array}
\label{eq:relq}
\end{equation}

Absolute $q$ sequence:
\begin{equation}
\begin{array}{l l l l l l}
                     \vdots & & \vdots  				& \vdots & \vdots & \\
j = \frac{f+1}{2} +k & & 0        				 &           0&          0&\cdots \\
j = \frac{f-1}{2} +k  & & H^0(B, H^{(f-1)/2 +k}(F)) & H^1(B, H^{(f-1)/2 +k}(F)) & H^2(B, H^{(f-1)/2 +k}(F)) 	
	&\cdots \\
\vdots && \vdots & \vdots & \vdots & \cdots \\
j =0 & & H^0(B, H^{0}(F)) & H^1(B, H^{0}(F)) & H^2(B, H^{0}(F)) &\cdots \\
E_2^{i,j}(\CY, q) && i=0&i=1&i=2& \cdots
\end{array}
\label{eq:absq}
\end{equation}

\noindent And using an analogous argument for the minimal complex, we get that the $E_2$ term of the 
relative $p$-perversity spectral sequence is:

\begin{equation}
\begin{array}{llllll}
j=f+1 & &H^0(B, H^{f}(F)) & H^1(B, H^{f}(F)) & H^1(B, H^{f}(F)) 	&\cdots \\
\vdots && \vdots                 & \vdots                  & \vdots                    & \cdots \\
j = \frac{f+3}{2} -k & & H^0(B, H^{(f+1)/2 -k}(F)) & H^1(B, H^{(f+1)/2 -k}(F)) & H^2(B, H^{(f+1)/2 -k}(F)) 	
	&\cdots \\
j = \frac{f+1}{2} -k & & 0&0&0&\cdots \\
\vdots && \vdots & \vdots & \vdots & \cdots \\
j =0 & & 0&0&0&\cdots \\
E_2^{i,j}(\CY, Y, p) && i=0&i=1&i=2& \cdots
\end{array}
\label{eq:relp}
\end{equation}

In the case that $f$ is even, these diagrams are similar, but with zeros
below level $f/2 + k$ for the relative $q$ spectral sequence, and shifted up by one as in the odd case,
zeros below $f/2-1- k$ for the relative $p$ spectral sequence, and again shifted up by one,
and zeros above level $f/2 + k$ for the absolute $q$ spectral sequence.

These are all truncated (and in the relative case, shifted) copies of the Leray-Serre spectral sequence for 
the fibration on $Y$.  So the higher levels of these spectral sequences will be the same as the 
higher levels of the spectral sequence for $Y$, but with term added because of the truncation.
By definition, $d_s \circ d_{s+1} = 0$ for any sequential differentials in the spectral sequence on $Y$.  Thus in the case of $E_{\infty}(I\!H_q^*(\CY))$ and $E_{\infty}(I\!H_p^*(\CY,Y)$, the extra terms will always have the form of $Im(d_s^{i,j})$ and will come from cases where there is no longer any image to quotient by because of the truncation above.  In the 
case of $E_{\infty}(I\!H_q^*(C_{\phi}Y,Y)$, the extra terms will have the form of a preimage $Im(d_s^{i,j})^*$ and will come from the 
cases where because of the truncation below, the kernel of the differential in the truncated sequence contains what would have been the image of the differential in the spectral sequence for $Y$.

It is easiest to understand this in an example.  Take the case when the base dimension $b=6$, the fibre dimension $f=5$ and the perversities are $p=\underline{\frakm} + 1$ and $q = \overline{\frakm} -1$.
Refer to the figures at the end of this paper for this example.
The the $E_{\infty}$ term of the spectral sequence for $I\!H_q^*(\CY)$ looks like figure \ref{fig:qabs}, 
the $E_{\infty}$ term of the spectral sequence for $I\!H_q^*(\CY, Y)$ looks like figure \ref{fig:qrel}
and the $E_{\infty}$ term of the spectral sequence for $I\!H_p^*(\CY, Y)$ looks like figure \ref{fig:prel}.
So in this example, we get for instance:
$$
E^{2,5}_\infty(I\!H_{\overline{\frakm} -1}^*(\CY,Y)) = E_{\infty}^{2,4}(Y) + Im(d_2^{2,4})^* + Im(d_3^{2,4})^*+
Im(d_4^{2,4})^*
$$
$$
E^{3,3}_\infty(I\!H_{\overline{\frakm} -1}^*(\CY)) = E_{\infty}^{3,3}(Y) + Im(d_2^{1,4}) + Im(d_3^{0,5})
$$
\noindent and 
$$
E^{2,5}_\infty(I\!H_{\underline{\frakm} + 1}^*(\CY,Y)) = E_{\infty}^{2,4} (Y)+ Im(d_4^{2,4})^*
$$
In general we get that:
$$
E^{i,j}_\infty(I\!H_q^*(\CY)) \cong
\left\{
\begin{array}{ll}
0 & j> f/2 +k \\
E_{\infty}^{ij}(Y) + \sum_{s \geq 1} Im(d_s^{i-s,j+s-1}) & j \leq f/2 + k, \\
&i+j > (f+1)/2 +k \\
E_{\infty}^{ij}(Y) & \mbox{ otherwise, }
\end{array} \right.
$$
$$
E^{i,j}_\infty(I\!H_q^*(\CY,Y)) \cong 
\left\{
\begin{array}{ll}
0 & j<(f+3)/2 +k \\
E_{\infty}^{i,j-1}(Y) + \sum_{s \geq j - f/2 - k} Im(d_s^{i,j-1})^* & j \geq (f+3)/2 + k,\\
& i+j < (f+1)/2+b +k \\
E_{\infty}^{i,j-1}(Y) & \mbox{ otherwise, }
\end{array} \right.
$$
\noindent and 
$$
E^{i,j}_\infty(I\!H_p^*(\CY,Y)) \cong 
\left\{
\begin{array}{ll}
0 & j<f/2 +1-k \\
E_{\infty}^{i,j-1}(Y) + \sum_{s \geq j - (f-1)/2 + k} Im(d_s^{i,j-1})^* & j \geq f/2 + 1-k, \\
&i+j < (f+1)/2+b -k \\
E_{\infty}^{i,j-1}(Y) & \mbox{ otherwise.}
\end{array} \right.
$$
When we map $I\!H_p^*(\CY,Y)$ to $I\!H_q^*(\CY,Y)$ we simply truncate the bottom of 
$E^{i,j}_\infty(I\!H_p^*(\CY,Y)) $ by $2k$ levels, if $f$ is odd, and by $2k+1$ levels if $f$
is even.  Modfying slightly the argument from \cite{CD} we have the 
\begin{proposition}
In terms of the spectral sequence, the map $I\!H_q^*(\CY,Y) \rightarrow I\!H_q^*(\CY)$
is zero on the terms of the form $E_{\infty}^{i,j}(Y)$ and is given by applying the appropriate
$d_r$'s to the other factors.  
\label{prop:incmap}
\end{proposition}
\begin{proof}
First consider a form $[\theta]$ in $Im(d_s^{i,j-1})^*$.  Then $d_s \theta$ will be in 
$E_s^{i+s, j-s}$ where $j-s$ is below the level of truncation for the spectral sequence.
That is, $j-s \leq (f-1)/2 + k$ if $f$ is odd and $f \leq f/2 + k$ if $f$ is even.  So $d \theta$ will have 
fibre degree $\leq f/2 + k$.  From Theorem \ref{th:hodge}, we know that in terms of 
maximal cohomology on $(\CY, g_c)$, $k = [[1+1/2c]]$ if $f$ is even and $k = [[1/2 + 1/2c]]$
if $f$ is odd.  So in either case, the fibre degree of $d\theta < f/2 + 1/2c$.  By Proposition 2, 
this guarantees that $d\theta$ extends to an $L^2$ form on $\CY$.  
In terms of the complex for relative maximal cohomology, the class
$[\theta, 0] \in \Omega^{*-1}(Y) \oplus \Omega^*(\CY, g_c)$ is equivalent to 
$[\theta, \alpha]$ in $E_s$ for any $\alpha in  \Omega^*(\CY, g_c)$, so choose 
$\alpha = d\theta$ to get a representative $[\theta, d\theta]$ which is in $I\!H_q^*(\CY,Y)$.
Then under the map to absolute cohomology, this goes to $[d\theta] \in I\!H_q^*(\CY)$.
The map from $\Omega^*(\CY,g_c)$ to $\Omega^*(Y)$ which induced the identification
of $I\!H_q^*(\CY)$ with spectral sequence terms for $Y$ is given by restriction.
Since $d\theta$ is constant in the $r$ direction, this is just $d_Y\theta$, which by 
identification with spectral sequence terms is just $d_s \theta$.

Now suppose that $[\theta] \in E_{\infty}^{i,j-1}(Y)$.  Then $d_Y \theta = 0$, so by the same 
argument as above, we get that $[\theta]$ goes to $[0]$ under the map
$I\!H_q^*(\CY,Y) \rightarrow I\!H_q^*(\CY)$.

\end{proof}

So we get that 
$$
Im(E_{\infty}^{ij}(\CY, Y, p) \rightarrow E_{\infty}^{ij}(\CY, q)) \cong
\left\{
\begin{array}{ll}
\sum_{s \geq j - (f-1)/2+k} Im(d_s^{i,j-1}) & i+j < (f+1)/2 +b-k ,\\
&j \geq (f+3)/2 +k\\
0 & \mbox{ otherwise}
\end{array}
\right. .
$$
\noindent
We obtain the result of the lemma by summing over $i+j=n/2=(b+f+1)/2$.
\end{proof}

Now we have to understand the signature pairing on
$Im(I\!H_p^{n/2}(C_{\phi}Y, Y) \rightarrow I\!H_q^{n/2}(C_{\phi}Y))$
in terms of this decomposition.  To do this, first return to our example where $f=5$, $b=6$ and $k=1$. 
We have
$$
I\!H_{\overline{\frakm}-1}^6(C_{\phi}Y)) \cong 
 E_{\infty}^{6,0}(I\!H^*_q(\CY))+ E_{\infty}^{5,1}(I\!H^*_q(\CY))+ E_{\infty}^{4,2}(I\!H^*_q(\CY))
+ E_{\infty}^{3,3}(I\!H^*_q(\CY))
$$
and
$$
Im(I\!H_{\underline{\frakm}+1}^6(C_{\phi}Y, Y) \rightarrow I\!H_{\overline{\frakm}-1}^6(C_{\phi}Y,Y))
$$
$$
\cong (Im(E_{\infty}^{0,6}(I\!H^*_p(\CY, Y)) \rightarrow E_{\infty}^{0,6}(I\!H^*_q(\CY,Y)))
+( Im(E_{\infty}^{1,5}(I\!H^*_p(\CY, Y)) \rightarrow E_{\infty}^{1,5}(I\!H^*_q(\CY,Y)) )
$$
Only forms with complementary bi-degrees for the the cone bundle $\CY$ can have nontrivial pairing.  By the filtration on Leray Serre spectral sequences, (see \cite{CD}), we can 
always represent a class in $E_s^{i,j}(\CY, q)$ or  $E_s^{i,j}(\CY, Y, p)$ by a form $\theta$ which
is a sum of forms of bidegree $(i-a, j+a)$ for $a >0$ and for which $d \theta$ is a sum of forms whose 
base degree is at least $j+s$.   Thus terms in the decomposition \ref{eq:decomp} whose $i$ indices add 
to more than the dimension of the base will have trivial intersection.  
Therefore we can represent the intersection form by a 4 by 2 block lower-triangular matrix:
$$
\left(
\begin{array}{llll}
A & 0 & 0 & 0\\
* & B & 0 & 0
\end{array}
\right),
$$
\noindent
where $A$ is the pairing between 
$Im(E_{\infty}^{0,6}(I\!H^*_p(\CY, Y)) \rightarrow E_{\infty}^{0,6}(I\!H^*_q(\CY,Y))$ and 
$E_{\infty}^{6,0}(I\!H^*_q(\CY))$ and $B$ is the pairing between 
$Im(E_{\infty}^{1,5}(I\!H^*_p(\CY, Y)) \rightarrow E_{\infty}^{1,5}(I\!H^*_q(\CY,Y))$ and 
$E_{\infty}^{5,1}(I\!H^*_q(\CY))$.

Now we have to understand these pieces.  First consider $A$.
We can further decompose $A$ into a 3 by 4 block matrix corresponding to the 
decompositions:
$$
E_{\infty}^{6,0}(I\!H^*_q(\CY)) \cong E_{\infty}^{6,0}(Y) + Im(d_6^{0,5})+ Im(d_5^{1,4}) .
$$
$$
(Im(E_{\infty}^{0,6}(I\!H^*_p(\CY, Y)) \rightarrow E_{\infty}^{0,6}(I\!H^*_q(\CY,Y)))
\cong 
E_{\infty}^{0,5}(Y) + Im(d_6^{0,5})^*+ Im(d_5^{0,5})^* + Im(d_4^{0,5})^* 
$$
We know that the entire intersection pairing descends to a nondegenerate pairing on 
$$
Im(I\!H_p^{n/2}(C_{\phi}Y, Y) \rightarrow I\!H_q^{n/2}(C_{\phi}Y))
$$
\noindent
which contains no terms of the form $E_{\infty}^{i,j}(Y)$ by Proposition \ref{prop:incmap}.
Thus such terms must pair trivially with everything.  So the first row and the first
column of block matrices in $A$ are trivial.

From the exact sequences in \ref{eq:fibreexactsequences}, we know that under the identifications
above, an element $[\theta] \in E_r^{i,j}(Y)$ lifts to the element 
$[d\theta, \theta] \in I\!H_{\underline{\frakm}+1}^6(C_{\phi}Y, Y)$, which is constant in the $r$-direction
on $\CY$, and an element $[\alpha] \in E_\infty^{i,j}(M) + \sum Im(d_r)$ lifts to 
$[\alpha] \in I\!H_q^{n/2}(C_{\phi}Y)$.  
Then for the lifts $[d\theta, \theta] \in I\!H_p^{n/2}(C_{\phi}Y, Y)$ and $\alpha \in I\!H_q^{n/2}(C_{\phi}Y)$, 
the pairing is given by 
$$
\langle [\theta],[\alpha] \rangle = \int_{\CY} d\theta \wedge \alpha - \int_M \theta \wedge \alpha.
$$
Neither $d\theta$ nor $\alpha$ contains a $dr$ term, so the first integral always vanishes.
Thus we are left with just the second term.

If the level of $[\theta]$ is greater than the level of $[\alpha]$ in the spectral sequence for $Y$, e.g. $[\theta] \in Im(d_6^{0,5})^*$ and $[\alpha] \in Im(d_5^{1,4})$, 
then we can choose $\theta$ to be a sum of forms with base degree greater than or equal to 0 and such that $d\theta$ has base degree greater than or equal to $i+r=6$ and fibre degree less than or equal to $j-r+1=0$. Similarly, we can choose $\alpha$ to have base degree greater than or equal to 6 and to equal $d_5\eta$ for a form of base degree greater than or equal to $1$ and fibre degree less than or equal to $4$.  So we have
$$
\langle [\theta],[\alpha] \rangle = -\int_Y \theta \wedge \alpha = -\int_Y \theta \wedge d_5 \eta
$$
\noindent Only the lowest base degree
parts of these forms will pair nondegenerately on $Y$, so we can add the other terms of $d \eta$
(which will all have higher base degree) to $d_5 \eta$ without changing the integral.  So we have
$$
= - \int_Y \theta \wedge d\eta = \int_Y d\theta \wedge \eta 
$$ 
and similarly, we can eliminate all but $d_6 \theta$ without changing the integral to get:
$$
= \int_Y d_6 \theta \wedge \eta.
$$

Now $d_6 \theta$ has base degree greater than or equal to 6
and $\eta$ has base degree greater than or equal to 1.  
So none of their component forms have complementary bidegree since the base degrees add to 
more than 6.  So this integral vanishes.  
Thus the block matrix for $A$ is also of lower triangular form:
$$
A = \left(
\begin{array}{lll}
0 & 0 & 0  \\
0 & A_1& 0 \\
0 & * & A_2 \\
0 & * & *
\end{array}
\right),
$$
\noindent 
where $A_1$ is the intersection matrix for $Im(d_6^{0,5})^* \otimes Im(d_6^{0,5})$
and $A_2$ is the intersection matrix for $Im(d_5^{0,5})^* \otimes Im(d_5^{1,4})$ .
The pairing $A_1$ is given by 
$$
Im(d_6^{0,5})^* \otimes Im(d_6^{0,5})^* \rightarrow \mathbb{R}
$$
$$
\omega  , \alpha \rightarrow -\langle \omega \cdot d_6 \alpha, \zeta_6 \rangle
$$
where $\zeta_6$ is the volume form in $E_6$.  
The signature of $A_1$ is therefore exactly $-\tau_6$.  

By a similar decomposition, $B$ has the form:
$$
B = \left(
\begin{array}{lll}
0 & 0& 0 \\
0 & B_1 & 0 \\
0 & * & B_2
\end{array}
\right),
$$
\noindent where $B_1$ is the intersection matrix for $Im(d_5^{1,4})^* \otimes  Im(d_5^{0,5})$
and $B_2$ is the intersection matrix for $Im(d_4^{1,4})^* \otimes  Im(d_4^{1,4})$.
The signature of $B_2$ is exactly $-\tau_4$, and the matrix $B_1 = - A_2$ since for 
any $\omega \in Im(d_5^{1,4})^*$ and $\alpha \in Im(d_5^{0,5})$, we have
$$
\langle [\omega],[\alpha] \rangle = - \int_Y \omega \wedge \alpha 
= - \int_Y \omega \wedge d_5 \eta 
$$
$$
= - \int_Y \omega \wedge d \eta = \int_Y d\omega \wedge \eta = \int_Y d_5 \omega \wedge \eta
$$
$$
= - \langle [\gamma], [\eta] \rangle
$$
for some classes $[\gamma] \in Im(d_5^{1,4})$ and $\eta \in Im(d_5^{0,5})^*$.
So the signatures of $B_1$ and $A_2$ cancel.  We can observe from this argument also that
in this case $\tau_5 = 0$, and so we are left with 
$$
\hat{\sigma}_p(\CY) = -\tau_4 - \tau_6 = \sum_{s=2+2k} \tau_s.
$$

The general argument is similar.  That is, since only complementary bidegrees will pair, we can always represent the intersection pairing on 
$Im(I\!H_p^{n/2}(C_{\phi}Y, Y) \rightarrow I\!H_q^{n/2}(C_{\phi}Y))$
by a block lower triangular matrix.  By integration by parts and the degeneracy of terms of the form
$E_{\infty}^{i,j}$ in the pairing, each block will further decompose into a block lower triangular matrix.  One of these sublocks will have signature $\tau_s$ for some $s$ 
and the rest will have signature which cancels the signature of one of the sub-blocks of 
another block in the large matrix.  The signature of the large matrix will thus be the 
sum of $\tau_s$'s.  If $f$ is odd, only even $\tau_s$ will contribute and if $f$ is even, only 
odd $\tau_s$ will contribute.  The first $s$ to appear comes from the lowest appropriate parity
$d_s$ which appears in 
$Im(E_{\infty}^{i,j} (I\!H_p^*(\CY, Y)) \rightarrow E_{\infty}^{i,j} (I\!H_q^*(\CY, Y))$,
which will be $s = 2 + 2k$ if $f$ is odd and $s = 3 + 2k$ if $f$ is even.  Since the opposite parity
$\tau_s$ vanish, we can simply write
$$
\hat{\sigma}_p(\CY) = - \sum_{s=2+ 2k} \tau_s
$$
and we are done.

\end{proof}


\section{Further work}
It seems likely that the theorems in this paper should generalize.  For instance, by relying more on strictly topological methods and using the basic definitions and results about intersection cohomology, it should be possible to define a signature pairing for any perversity intersection cohomology on any pseudomanifold and to prove a more general version of Novikov additivity:
\begin{conjecture}
If $X^n$ is a pseudomanifold and if $Y \subset X$ is a compact codimension 1 submanifold such that 
$X = Z \bigcup_Y Z^{\prime}$ where
$Y \subset \subset X^{\mbox{reg}}$, then for any perversity function, $\mathfrak{p}$, the signature
of the intersection pairing defined on $IH^{n/2}_{\mathfrak{p}}(X)$ satisfies:
$$
\sigma_{\mathfrak{p}}(X) = \hat{\sigma}_{\mathfrak{p}}(Z) + \hat{\sigma}_{\mathfrak{p}}(Z^{\prime}),
$$
\noindent
where as before, $\hat{\sigma}_{\mathfrak{p}}(Z)$ is the signature for the cohomology that satisfies
perversity $\mathfrak{p}$ type conditions away from $Y$ and relative boundary conditions at $Y$.
\end{conjecture}

Together with Daniel Grieser the author is currently proving the following conjecture which also relates the topological signatures considered in this paper to $L^2$ signatures for a family of complete metrics on $M$ that interpolates between fibred cusp and cylindrical metrics:

\begin{conjecture} If $M$ is a manifold with boundary fibration 
$\del M \stackrel{\phi}{\rightarrow} B$ 
and the metric on
$M$ is quasi-isometric near the boundary to one of the form
$$
ds_c^2 = R^{2c}(dR^2 + k) + \phi^*ds_B^2,
$$
where $R \in [1, \infty)$, $-1 \leq c<0$ and $k$ is a symmetric two-tensor on $\del M$ which restricts
to a metric on each fiber, then for $m$ even, 
$$
\mathcal{H}^i_{(2)}(M,g_c) \cong 
 \mbox{Im} \, \big(I\!H^k_{ \underline{\frakm} + [[-1/2c]]}(X,\hB)
 \longrightarrow I\!H^k_{ \overline{\frakm}  - [[-1/2c]]}(X,\hB) \big)
 $$
 and for $m$ odd, 
 $$
\mathcal{H}^i_{(2)}(M,g_c) \cong 
\mbox{Im} \, \big(I\!H^k_{ \underline{\frakm}  +[[(-1/2c)-(1/2)]]}(X,\hB)
\longrightarrow I\!H^k_{ \overline{\frakm}  -[[(-1/2c)-(1/2)]]}(X,\hB) \big),
$$
\noindent where $[[x]]$ denotes the smallest integer strictly greater than $x$.
Thus we obtain the signature theorem:
$$
L^2-\sigma(M)= \mbox{sgn}\,
\mbox{Im}\,\big( H^*(M, \del M) \longrightarrow H^*(M))\big) + \sum_{i=2+2k}^{\infty} \tau_i,
$$
\noindent
where $k= [-1/2c]$ when $m$ is even and $[(-1/2c)-(1/2)]$ when $m$ is odd.

\end{conjecture}

It seems likely that $IH^*_{\mathfrak{p}}(X)$ should be isomorphic to the space of $L^2$-harmonic forms on $X^{reg}$ for an incomplete metric on $X^{reg}$ with degeneration conditions near the singular strata related to those for the family of metrics considered in this paper.  It would be interesting to explore 
which metrics these might be, and also to consider if it might be possible to define a general 
$\tau$-invariant for singular sets of pseudomanifolds with more than one singular stratum.  Such a 
generalization would make it possible to derive signature theorems from the Hodge theorems
for noncompact manifolds with more complicated singularity structures than we have so far considered.

\end{document}